\documentclass[12pt,reqno]{amsart}
\usepackage{amsmath,amssymb}
\usepackage[dvips]{graphics}
\usepackage{epsfig}

\newcommand{\ncr}[2]{{#1 \choose #2}}
\newtheorem{thm}{Theorem}[section]

\newtheorem{lem}[thm]{Lemma}

\newtheorem{exa}[thm]{Example}
\newtheorem{defi}[thm]{Definition}

\numberwithin{equation}{section}

\newtheorem{rek}[thm]{Remark}
\theoremstyle{definition}

\newcommand\ben{\begin{enumerate}}
\newcommand\een{\end{enumerate}}

\newcommand{\twocase}[5]{#1 \begin{cases} #2 & \text{#3}\\ #4
&\text{#5} \end{cases}   }

\newcommand\be{\begin{equation}}
\newcommand\ee{\end{equation}}
\newcommand\bea{\begin{eqnarray}}
\newcommand\eea{\end{eqnarray}}

\newcommand{\R}{\mathbb{R}}

\newcommand{\Z}{\mathbb{Z}}

\newcommand{\gep}{\epsilon}  %lowercase epsilon
\newcommand{\E}{{\mathbb E}} % expectation
\numberwithin{equation}{section}

\begin{document}

\title[The Modulo $1$ CLT and Benford's Law for Products]{The
Modulo $1$ Central Limit Theorem and Benford's Law for Products}
\author{Steven J. Miller}
\email{sjmiller@math.brown.edu}\address{Department of Mathematics,
Brown University, Providence, RI $02912$}
\thanks{We thank Brian Cole, Ted Hill, Jill Pipher and Sergei Treil for discussions
on an earlier draft. The first author was partially supported by NSF
grant DMS-0600848.}

\author{Mark
J. Nigrini} \email{mnigrini@smcvt.edu}
\address{Department of Business Administration and Accounting, Saint Michael's College,
Colchester, VT 05439} \subjclass[2000]{60F05, 60F25, 11K06
(primary), 42A10, 42A61, 62E15 (secondary).}\keywords{Central Limit
Theorem, sums modulo $1$, Fourier series, $L^1$-convergence,
Benford's Law}
% 11B83
\date{\today}

%and let $X_{1:N}, \dots, X_{N:N}$ be the $X_m$'s ranked in
%increasing order.

%If $B=e$ this distribution is \emph{almost} the same as Benford's
%Law base $e$; the difference from Benford behavior is more
%pronounced base $10$, though still quite small.

\begin{abstract} Using elementary results from Fourier analysis,
we provide an alternate proof of a necessary and sufficient
condition for the sum of $M$ independent continuous random variables
modulo $1$ to converge to the uniform distribution in $L^1([0,1])$,
and discuss generalizations to discrete random variables. A
consequence is that if $X_1, \dots, X_M$ are independent continuous
random variables with densities $f_1, \dots, f_M$, for any base $B$
as $M \to \infty$ for many choices of the densities the distribution
of the digits of $X_1 \cdots X_M$ converges to Benford's law base
$B$. The rate of convergence can be quantified in terms of the
Fourier coefficients of the densities, and provides an explanation
for the prevalence of Benford behavior in many diverse systems. To
highlight the difference in behavior between identically and
non-identically distributed random variables, we construct a
sequence of densities $\{f_i\}$ with the following properties: (1)
for each $i$, if every $X_k$ is independently chosen with density
$f_i$ then the sum converges to the uniform distribution; (2) if the
$X_k$'s are independent but non-identical, with $X_k$ having
distribution $f_k$, then the sum does not converge to the uniform
distribution.
\end{abstract}

\maketitle

%%%%%%%%%%%%%%%%%%%%%%%%%%%%%%%%%%%%%%%%%%%%%%%%%%%%%%%%%%%%%%%%%%%%%%%%%%%%%%%%%%
%%%%%%%%%%%%%%%%%%%%%%%%%%%%%%%%%%%%%%%%%%%%%%%%%%%%%%%%%%%%%%%%%%%%%%%%%%%%%%%%%%

\section{Introduction}

We investigate necessary and sufficient conditions for the
distribution of a sum of random variables modulo $1$ to converge to
the uniform distribution. This topic has been fruitfully studied by
many previous researchers. Our purpose here is to provide an
elementary proof of prior results, and explicitly connect this
problem to related problems in the Benford's Law literature
concerning the distribution of the leading digits of products of
random variables. As this question has motivated much of the
research on this topic, we briefly describe that problem and its
history, and then state our results.

For any base $B$ we may uniquely write a positive $x\in\R$ as $x =
M_B(x)\cdot B^k$, where $k\in \Z$ and $M_B(x)$ (called the mantissa)
is in $[1,B)$. A sequence of positive numbers $\{a_n\}$ is said to
be \textbf{Benford base $B$} (or to satisfy Benford's Law base $B$)
if the probability of observing the base-$B$ mantissa of $a_n$ of at
most $s$ is $\log_B s$. More precisely, \be \lim_{N \to \infty}
\frac{ \#\{n \le N: \text{$1 \le M_B(a_n) \le s$} \} }{N} \ =\
\log_B s. \ee Benford behavior for continuous systems is defined
analogously. Thus base $10$ the probability of observing a first
digit of $j$ is $\log_{10} (j+1) - \log_{10} (j)$, implying that
about $30\%$ of the time the first digit is a $1$.

Benford's Law was first observed by Newcomb in the 1880s, who
noticed that pages of numbers starting with $1$ in logarithm tables
were significantly more worn than those starting with $9$. In 1938
Benford \cite{Ben} observed the same digit bias in 20 different
lists with over 20,000 numbers in all. See \cite{Hi1,Rai} for a
description and history. Many diverse systems have been shown to
satisfy Benford's law, ranging from recurrence relations \cite{BrDu}
to $n!$ and $\ncr{n}{k}$ ($0 \le k \le n$) \cite{Dia} to iterates of
power, exponential and rational maps \cite{BBH,Hi2} to values of
$L$-functions near the critical line and characteristic polynomials
of random matrix ensembles \cite{KoMi} to iterates of the $3x+1$ Map
\cite{KoMi,LS} to differences of order statistics \cite{MN}. There
are numerous applications of Benford's Law. It is observed in
natural systems ranging from hydrology data \cite{NM} to stock
prices \cite{Ley}, and is used in computer science in analyzing
round-off errors (see page 255 of \cite{Knu} and \cite{BH}), in
determining the optimal way to store numbers\footnote{If the data is
distributed according to Benford's Law base $2$, the probability of
having to shift the result of multiplying two numbers if the
mantissas are written as $0.x_1x_2x_3\cdots$ is about $.38$; if they
are written as $x_1.x_2x_3\cdots$ the probability is about $.62$.}
\cite{Ha}, and in accounting to detect tax fraud \cite{Nig1,Nig2}.
See \cite{Hu} for a detailed bibliography of the field.

In this paper we consider the distribution of digits of products of
independent random variables, $X_1 \cdots X_M$, and the related
questions about probability densities of random variables modulo
$1$. Many authors \cite{Sa,ST,AS,Adh,Ha,Tu} have observed that the
product (and more generally, any nice arithmetic operation) of two
random variables is often closer to satisfying Benford's law than
the input random variables; further, that as the number of terms
increases, the resulting expression seems to approach Benford's Law.

Many of the previous works are concerned with determining exact
formulas for the distribution of $X_1 \cdots X_M$; however, to
understand the distribution of the digits all we need is to
understand $\log_B |X_1 \cdots X_M| \bmod 1$. This leads to the
equivalent problem of studying sums of random variables modulo $1$.
This formulation is now ideally suited for Fourier analysis. The
main result is a variant of the Central Limit Theorem, which in this
context states that for ``nice'' random variables, as $M\to\infty$
the sum\footnote{That is, we study sums of the form $Y_1 + \cdots +
Y_M$. For the standard Central Limit Theorem one studies
$\frac{\sum_m Y_m - \E[\sum_m Y_m]}{{\rm StDev}(\sum_m Y_m)}$. We
subtract the mean and divide by the standard deviation to obtain a
quantity which will be finite as $M\to\infty$; however, sums modulo
$1$ are a priori finite, and thus their unscaled value is of
interest.} of $M$ independent random variables modulo $1$ tends to
the uniform distribution; by simple exponentiation this is
equivalent to Benford's Law for the product (see \cite{Dia}). To
emphasize the similarity to the standard Central Limit Theorem and
the fact that our sums are modulo $1$, we refer to such results as
Modulo $1$ Central Limit Theorems. Many authors
\cite{Bh,Bo,Ho,JR,Lev,Lo,Ro,Sc1,Sc2,Sc3} have analyzed this problem
in various settings and generalizations, obtaining sufficient
conditions on the random variables (often identically distributed)
as well as estimates on the rate of convergence.

Our main result is a proof, using only elementary results from
Fourier analysis, of a necessary and sufficient condition for a sum
modulo 1 to converge to the uniform distribution in $L^1([0,1])$. We
also give a specific example to emphasize the different behavior
possible when the random variables are not identically distributed.
We let $\widehat{g_{m}}(n)$ denote the $n$\textsuperscript{th}
Fourier coefficient of a probability density $g_m$ on $[0,1]$: \be
\widehat{g_{m}}(n) \ = \ \int_0^1 g_m(x) e^{-2\pi i n x} dx.\ee

\begin{thm}[The Modulo $1$ Central Limit Theorem for Independent
Continuous Random Variables]\label{thm:mainL1} Let $\{Y_m\}$ be
independent continuous random variables on $[0,1)$, not necessarily
identically distributed, with densities $\{g_m\}$. A necessary and
sufficient condition for the sum $Y_1 + \cdots + Y_M$ modulo $1$ to
converge to the uniform distribution as $M\to\infty$ in $L^1([0,1])$
is that for each $n\neq 0$ we have $\lim_{M\to\infty}
\widehat{g_1}(n)\cdots \widehat{g_M}(n) = 0$.
\end{thm}

As Benford's Law is equivalent to the associated base $B$ logarithm
being equidistributed modulo $1$ (see \cite{Dia}), from Theorem
\ref{thm:mainL1} we immediately obtain the following result on the
distribution of digits of a product.

\begin{thm}\label{thm:mainbenfL1} Let $X_1, \dots, X_M$ be independent
continuous random variables, and let $g_{B,m}$ be the density of
$\log_B M_B(|X_m|)$. A necessary and sufficient condition for the
distribution of the digits of $X_1 \cdots X_M$ to converge to
Benford's Law (base $B$) as $M\to\infty$ in $L^1([0,1])$ is for each
$n\neq 0$ that $\lim_{M\to\infty} \widehat{g_{B,1}}(n)\cdots
\widehat{g_{B,M}}(n) = 0$. \end{thm}

As other authors have noticed, the importance of results such as
Theorem \ref{thm:mainbenfL1} is that they give an explanation of why
so many data sets follow Benford's Law (or at least a close
approximation to it). Specifically, if we can consider the observed
values of a system to be the product of many independent processes
with reasonable densities, then the distribution of the digits of
the resulting product will be close to Benford's Law.

We briefly compare our approach with other proofs of results such as
Theorem \ref{thm:mainL1} (where the random variables are often taken
as identically distributed). If the random variables are identically
distributed with density $g$, our condition reduces to
$|\widehat{g}(n)| < 1$ for $n \neq 0$. For a probability
distribution, $|\widehat{g}(n)| = 1$ for $n \neq 0$ if and only if
there exists $\alpha \in R$ such that all the mass is contained in
the set $\{\alpha, \alpha + \frac1n, \dots, \alpha + \frac{n-1}n\}$.
(As we are assuming our random variables are continuous and not
discrete, the corresponding densities are in $L^1([0,1])$ and this
condition is not met; in Theorem \ref{thm:discmainL1} we discuss
generalizations to discrete random variables.) In other words, the
sum of identically distributed random variables modulo $1$ converges
to the uniform distribution if and only if the support of the
distribution is not contained in a coset of a finite subgroup of the
circle group $[0,1)$. Interestingly, Levy \cite{Lev} proved this
just one year after Benford's paper \cite{Ben}, though his paper
does not study digits. Levy's result has been generalized to other
compact groups, with estimates on the rate of convergence \cite{Bh}.
Stromberg \cite{Str} proved that\footnote{The following formulation
is taken almost verbatim from the first paragraph of \cite{Bh}.}
\emph{the $n$-fold convolution of a regular probability measure on a
compact Hausdorff group $G$ converges to the normalized Haar measure
in the weak-star topology if and only if the support of the
distribution is not contained in a coset of a proper normal closed
subgroup of $G$.}

Our arguments in the proof of Theorem \ref{thm:mainL1} may be
generalized to independent discrete random variables, at the cost of
replacing $L^1$-convergence with weak convergence. Below
$\delta_{\alpha}(x)$ denotes a unit point mass at $\alpha$.

\begin{thm}[Modulo $1$ Central Limit Theorem for Certain Independent Discrete Random
Variables]\label{thm:discmainL1} Let $\{Y_m\}$ be independent
discrete random variables on $[0,1)$, not necessarily identically
distributed, with densities \be g_m(x) \ = \ \sum_{k=1}^{r_m}
w_{k,m} \delta_{\alpha_{k,m}}(x), \ \ \ w_{k,m}
> 0, \ \ \ \sum_{k=1}^{r_m} w_{k,m} \ = \ 1. \ee Assume that there is
a \emph{finite} set $A \subset [0,1)$ such that all $\alpha_{k,m}
\in A$. A necessary and sufficient condition for the sum $Y_1 +
\cdots + Y_M$ modulo $1$ to converge weakly to the uniform
distribution as $M\to\infty$ is that for each $n\neq 0$ we have
$\lim_{M\to\infty} \widehat{g_1}(n)\cdots \widehat{g_M}(n) = 0$.
\end{thm}

In \S\ref{sec:mainproof} we prove Theorem \ref{thm:mainL1} using
only elementary facts from Fourier analysis, showing our condition
is a consequence of Lebesgue's Theorem (on $L^1$-convergence of the
Fej\'{e}r series) and a standard approximation argument. We give an
example of distinct densities $\{f_i\}$ with the following
properties: (1) for each $i$, if every $X_k$ is independently chosen
with density $f_i$ then the sum converges to the uniform
distribution; (2) if the $X_k$'s are independent but non-identical,
with $X_k$ having distribution $f_k$, then the sum does not converge
to the uniform distribution. This example illustrates the difference
in behavior when the random variables are not identically
distributed: to obtain uniform behavior for the sum it does not
suffice for each random variable to satisfy Levy or Stromberg's
condition (the distribution is not concentrated on a coset of a
finite subgroup of $[0,1)$). We conclude in \S\ref{sec:discrete} by
sketching the proof of Theorem \ref{thm:discmainL1}, and in Appendix
\ref{sec:alternate} we comment on alternate techniques to prove
results such as Theorem \ref{thm:mainbenfL1} (in particular, why our
arguments are more general than applying the standard Central Limit
Theorem to $\log_B|X_1| + \cdots +\log_B|X_M|$ to analyze the
distribution of digits of $|X_1 \cdots X_N|$).

%We discuss in \S\ref{sec:alternate} by analyzing the limits on what
%the standard Central Limit Theorem can say on sums of random
%variables modulo $1$.

%%%%%%%%%%%%%%%%%%%%%%%%%%%%%%%%%%%%%%%%%%%%%%%%%%%%%%%%%%%%%%%%%%%%%%%%%%%%%%%%%%
%%%%%%%%%%%%%%%%%%%%%%%%%%%%%%%%%%%%%%%%%%%%%%%%%%%%%%%%%%%%%%%%%%%%%%%%%%%%%%%%%%

\section{Analysis of Sums of Continuous Random Variables}\label{sec:mainproof}

We recall some standard facts from Fourier analysis (see for example
\cite{SS}). The convolution of two functions in $L^1([0,1])$ is \be
(f \ast g)(x) \ = \ \int_0^1 f(y) g(x-y)dy \ = \ \int_0^1
f(x-y)g(y)dy. \ee Convolution is commutative and associative, and
the $n$\textsuperscript{th} Fourier coefficient of a convolution is
the product of the two $n$\textsuperscript{th} Fourier coefficients.

Let $g_1$ and $g_2$ be two probability densities in $L^1([0,1])$. If
$Z_i$ is a random variable on $[0,1)$ with density $g_i$, then the
density of $Z_1 + Z_2 \bmod 1$ is the convolution of $g_1$ with
$g_2$.

\begin{defi}[Fej\'{e}r kernel, Fej\'{e}r series] Let $f \in
L^1([0,1])$. The $N$\textsuperscript{{\rm th}} Fej\'{e}r kernel is
\be F_N(x) \ = \ \sum_{n=-N}^N \left(1 - \frac{|n|}{N}\right)
e^{2\pi i n x}, \ee and the $N$\textsuperscript{{\rm th}} Fej\'{e}r
series of $f$ is \be T_Nf(x) \ = \ (f \ast F_N)(x) \ = \
\sum_{n=-N}^N \left(1 - \frac{|n|}{N}\right) \widehat{f}(n) e^{2\pi
i n x}. \ee The Fej\'{e}r kernels are an approximation to the
identity (they are non-negative, integrate to $1$, and for any
$\delta \in (0,1/2)$ we have
$\lim_{N\to\infty}\int_{\delta}^{1-\delta} F_N(x)dx = 0$).\end{defi}

%\begin{thm}[Fej\'{e}r's Theorem]\label{thm:fejer}
%If $f$ is continuous and periodic
%on $[0,1]$ then as $N\to\infty$, $T_N f$ converges uniformly to $f$.
%\end{thm}

\begin{thm}[Lebesgue's Theorem]\label{thm:fejTnconverg}
Let $f\in L^1([0,1])$. As $N\to\infty$, $T_N f$ converges to $f$ in
$L^1([0,1])$. \end{thm}

\begin{lem}\label{lem:fejconvolve} Let $f, g \in L^1([0,1])$.
Then $T_N (f \ast g) = (T_N f) \ast g$. \end{lem}

\begin{proof} The proof follows immediately from the commutative and
associative properties of convolution. \end{proof}

%\begin{proof} As $\widehat{f \ast g}(n) =
%\widehat{f}(n)\widehat{g}(n)$, we have \bea (T_N (f \ast g))(x) \ =
%\ \sum_{n=-N}^N \left(1 - \frac{|n|}{N}\right) \widehat{f}(n)
%\widehat{g}(n) e^{2\pi i n x}. \eea Further \bea ((T_N f) \ast g)(x)
%& \ = \ & \int_0^1 \sum_{n=-N}^N \left(1 - \frac{|n|}{N}\right)
%\widehat{f}(n) e^{2\pi i n (x-y)} g(y) dy \nonumber\\ & = &
%\sum_{n=-N}^N \left(1 - \frac{|n|}{N}\right) \widehat{f}(n) e^{2\pi
%i n x} \int_0^1 g(y) e^{-2\pi i n y} dy \nonumber\\ & = &
%\sum_{n=-N}^N \left(1 - \frac{|n|}{N}\right) \widehat{f}(n)
%\widehat{g}(n) e^{2\pi i n x}; \eea the interchange of summation and
%integration is permissible as the sum is finite.
%\end{proof}

We can now prove Theorem \ref{thm:mainL1}.

\begin{proof}[Proof of Theorem \ref{thm:mainL1}] We first show our
condition is sufficient. The density of the sum modulo $1$ is $h_M =
g_1 \ast \cdots \ast g_M$. It suffices to show that, for any $\gep >
0$, \be \lim_{M\to\infty} \int_0^1 |h_M(x) - 1| dx \ < \ \gep. \ee
Using Lebesgue's Theorem (Theorem \ref{thm:fejTnconverg}), choose
$N$ sufficiently large so that \be\label{eq:h1Tn} \int_0^1 |h_1(x) -
T_N h_1(x)| dx \ < \ \frac{\gep}2. \ee

While $N$ was chosen so that \eqref{eq:h1Tn} holds with $h_1$, in
fact this $N$ works for \emph{all} $h_M$ (with the same $\gep$).
This follows by induction. The base case is immediate (this is just
our choice of $N$). Assume now that \eqref{eq:h1Tn} holds with $h_1$
replaced by $h_M$; we must show it holds with $h_1$ replaced by
$h_{M+1} = h_M \ast g_{M+1}$. By Lemma \ref{lem:fejconvolve} we have
\be T_N h_{M+1} \ = \ T_N (h_M \ast g_{M+1}) \ = \ (T_N h_M) \ast
g_{M+1}. \ee This implies \bea & & \int_0^1 |h_{M+1}(x) - T_N
h_{M+1}(x)| dx\nonumber\\ & &\ \ \ \ \ \ =\ \int_0^1 |(h_M \ast
g_{M+1})(x) - (T_N h_M)\ast g_{M+1}(x)| dx \nonumber\\ & &\ \ \ \ \
\ =\ \int_0^1 \left| \int_0^1 \left(h_M(y) - T_N h_M(y)\right) \cdot
g_{M+1}(x-y) \right| dy dx \nonumber\\ & &\ \ \ \ \ \ \le \ \int_0^1
\int_0^1 |h_M(y) - T_Nh_M(y)| \cdot g_{M+1}(x-y)dx dy \nonumber\\ &
&\ \ \ \ \ \ =\ \int_0^1|h_M(y) - T_Nh_M(y)|dy \cdot 1 \ < \
\frac{\gep}2; \eea the interchange of integration above is justified
by the absolute value being integrable in the product measure, and
the $x$-integral is $1$ as $g_{M+1}$ is a probability density.

To show $h_M$ converges to the uniform distribution in $L^1([0,1])$,
we must show $\lim_{M\to\infty}\int_0^1 |h_M(x) - 1|dx  = 0$. Let
$N$ and $\gep$ be as above. By the triangle inequality we have \bea
\int_0^1 |h_M(x) - 1|dx & \ \le \ & \int_0^1 |h_M(x) - T_N h_M(x)|dx
+ \int_0^1 |T_Nh_M(x) - 1|dx. \nonumber\\ \eea From our choices of
$N$ and $\gep$, $\int_0^1 |h_M(x) - T_N h_M(x)|dx < \gep/2$; thus we
need only show $\int_0^1 |T_Nh_M(x) - 1|dx < \gep/2$ to complete the
proof. As $\widehat{h_M}(0) = 1$, \bea \int_0^1 |T_Nh_M(x) - 1|dx &
\ = \ & \int_0^1 \left|\sum_{n=-N \atop n \neq 0}^N \left(1 -
\frac{|n|}{N}\right) \widehat{h_M}(n) e^{2\pi i n x}\right| dx
\nonumber\\ & \le & \sum_{n=-N \atop n \neq 0}^N \left(1 -
\frac{|n|}{N}\right) |\widehat{h_M}(n)|. \eea However,
$\widehat{h_M}(n) = \widehat{g_1}(n) \cdots \widehat{g_M}(n)$, and
by assumption tends to zero as $M\to\infty$ (as each
$\widehat{g_m}(n)$ is at most $1$ in absolute value, for each $n$
the absolute value of the product is non-increasing in $M$). For
fixed $N$ and $\gep$, we may choose $M$ sufficiently large so that
$|\widehat{h_M}(n)| < \gep/4N$ whenever $n \neq 0$ and $|n| \le N$.
Thus \be \int_0^1 |T_Nh_M(x) - 1|dx \ < \ 2N \cdot \frac{\gep}{4N} \
= \ \frac{\gep}{2}, \ee which implies \be \int_0^1 |h_M(x) - 1| dx \
< \ \gep\ee for $M$ sufficiently large. As $\gep$ is arbitrary, this
completes the proof of the sufficiency; we now prove this condition
is necessary.

Assume for some $n_0 \neq 0$ that $\lim_{M\to\infty}
|\widehat{h}_M(n_0)| \neq 0$ (where as always $h_M = g_1 \ast \cdots
g_M$). As the $g_m$ are probability densities, $|\widehat{g}_m(n)|
\le 1$; thus the sequence $\{|\widehat{h}_M(n)|\}_{M=1}^\infty$ is
non-increasing for each $n$, and hence by assumption converges to
some number $c_n \in (0,1]$.

Let $E_M(x) = h_M(x) - 1$; note $\widehat{E}_M(n) =
\widehat{h}_M(n)$ for $n\neq 0$. To show $h_M$ does not converge to
the uniform distribution on $[0,1]$, it suffices to show that $E_M$
does not converge almost everywhere to the zero function on $[0,1]$.
Let $n_0$ be as above. We have \bea \left|\widehat{h}_M(n_0)\right|
\ = \ \left|\widehat{E}_M(n_0)\right| \ = \ \left|\int_0^1 E_M(x)
e^{2\pi i n_0 x} dx \right| \ \ge \ c_{n_0} \ > \ 0. \eea Therefore
at least one of the following integrals is at least $c_{n_0}/2$:
\bea \int_{x \in [0,1] \atop {\rm Re}\left(E_M(x)\right) \ge 0} {\rm
Re}\left(E_M(x)\right)dx, \ \ & & \ \ \int_{x \in [0,1] \atop {\rm
Re}\left(E_M(x)\right) \le 0} {\rm
Re}\left(-E_M(x)\right)dx \nonumber\\
\int_{x \in [0,1] \atop {\rm Im}\left(E_M(x)\right) \ge 0} {\rm
Im}\left(E_M(x)\right)dx, \ \ & & \ \ \int_{x \in [0,1] \atop {\rm
Im}\left(E_M(x)\right) \le 0} {\rm Im}\left(-E_M(x)\right)dx, \eea
and $h_M$ cannot converge to the zero function in $L^1([0,1])$;
further, we obtain an estimate on the $L^1$-distance between the
uniform distribution and $h_M$.
\end{proof}

The behavior is non-Benford if the conditions of Theorem
\ref{thm:mainbenfL1} are violated. It is enough to show that we can
find a sequence of densities $g_{B,m}$ such that $\lim_{M\to\infty}
\prod_{m=1}^M\widehat{g_{B,m}}(1) \neq 0$. We are reduced to
searching for an infinite product that is non-zero; we also need
each term to be at most $1$, as the Fourier coefficients of a
probability density are dominated by $1$. A standard example is
$\prod_m c_m$, where $c_m = \frac{m^2+2m}{(m+1)^2}$; the limit of
this product is $1/2$. Thus as long as $\widehat{g_{B,m}}(1) \ge
\frac{m^2+2m}{(m+1)^2}$, the conclusion of Theorem
\ref{thm:mainbenfL1} will not hold for the products of the
associated random variables; analogous reasoning yields a sum of
independent random variables modulo $1$ which does not converge to
the uniform distribution.

\begin{exa}[Non-Benford Behavior of Products]\label{exa:nonbenf}
Consider \be \twocase{\phi_{m} \ = \ }{m}{{\rm if} $|x - \frac18|
\le \frac1{2m}$}{0}{{\rm otherwise;}} \ee $\phi_{m}$ is non-negative
and integrates to $1$. As $m \to \infty$ we have
$|\widehat{\phi_{m}}(1)| \to 1$ because the density becomes
concentrated at $1/8$ (direct calculation gives $\widehat{\phi_m}(1)
= e^{2\pi i/8} + O(m^{-2})$). Let $X_1, \dots, X_M$ be independent
random variables where the associated densities $g_{B,m}$ of $\log_B
M(|X_m|)$ are $\phi_{11^m}$. The behavior is non-Benford (see Figure
\ref{fig:counterex}). Note, however, that if each $X_m$ had the
common distribution $\phi_i$ for any fixed $i$, then in the limit
the product will satisfy Benford's law.
\begin{figure}
\begin{center}
\scalebox{.75}{\includegraphics{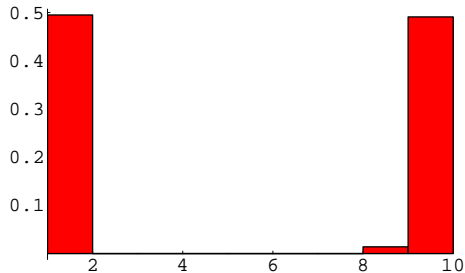}}
\caption{\label{fig:counterex} Distribution of digits (base 10) of
1000 products $X_1 \cdots X_{1000}$, where $g_{10,m} =
\phi_{11^m}$.}
\end{center}\end{figure}
\end{exa}

\begin{rek} Generalizations of Theorem \ref{thm:mainL1} hold for more general sums of random
variables. Instead of $Y_1 + \cdots + Y_M$ we may study $\eta_1 Y_1
+ \cdots + \eta_M Y_M$, where each $\eta_m$ is a random variable
taking values in $\{-1,1\}$; the proof follows from the observation
that if $Y_m$ has density $g_m(y)$ then $-Y_m$ has density
$g_m(1-y)$.
\end{rek}

%%%%%%%%%%%%%%%%%%%%%%%%%%%%%%%%%%%%%%%%%%%%%%%%%%%%%%%%%%%%%%%%%%%%%%%%%%%%%%%%%%
%%%%%%%%%%%%%%%%%%%%%%%%%%%%%%%%%%%%%%%%%%%%%%%%%%%%%%%%%%%%%%%%%%%%%%%%%%%%%%%%%%

\section{Analysis of Sums of Discrete Random Variables}\label{sec:discrete}

Many results from Fourier analysis do not apply if the random
variables are discrete; Lebesgue's Theorem cannot be correct for a
point mass as the density is concentrated on a set of measure zero.
Let $\delta_\alpha(x)$ be a unit point mass\footnote{Thus
$\delta_\alpha(x)$ is a Dirac delta functional; if $\phi(x)$ is a
Schwartz function then $\int_0^1 \delta_\alpha(x)\phi(x)dx$ is
defined to be $\phi(\alpha)$.} at $\alpha$. Its Fourier coefficients
are $\widehat{\delta_\alpha}(n) = e^{-2\pi i n \alpha}$, and simple
algebra shows that its Fej\'{e}r series is \be F_N \delta_\alpha(x)
\ = \ \frac{e^{-2\pi i(N-1)(x-\alpha)}(e^{2\pi i
N(x-\alpha)}-1)^2}{(e^{2\pi i(x-\alpha)}-1)^2\ N}.\ee For $x\neq
\alpha$, $\lim_{N\to\infty} F_N \delta_\alpha(x)$ $=$
$\delta_\alpha(x)$ $=$ $0$; moreover, for $x$ near $\alpha$ we have
$|F_N \delta_\alpha(x)| \sim N$. Instead of convergence in
$L^1([0,1])$ we have weak convergence: for any Schwartz function
$\phi$,  \be \lim_{N\to\infty} \int_0^1 F_N \delta_\alpha(x)
\phi(x)dx \ = \ \int_0^1 \delta_\alpha(x) \phi(x)dx \ = \
\phi(\alpha). \ee

\begin{proof}[Sketch of the proof of Theorem \ref{thm:discmainL1}]
We argue as in Theorem \ref{thm:mainL1}.
Note Lemma \ref{lem:fejconvolve} holds if $f$ and $g$ are sums of
point masses. Instead of using Lebesgue's Theorem, we use weak
convergence: given an $\epsilon>0$ and a Schwartz function
$\phi(x)$, by weak convergence there is an $N$ such that
\be\label{eq:genh1Tn} \left| \int_0^1 \left(h_1(x) -
T_Nh_1(x)\right) \phi(x)dx \right| \ < \ \frac{\gep}2. \ee This is
the generalization of \eqref{eq:h1Tn}. Further, we may assume
\eqref{eq:genh1Tn} holds with $\phi(x)$ replaced with
$\phi_{\alpha_{k,m}}(x) = \phi(x+\alpha_{k,m})$ for any
$\alpha_{k,m} \in A$. \emph{This is only true because $A$ is
finite}; while $N=N(\phi)$ depends on $\phi$, as there are  only
finitely many test functions $\phi_{\alpha_{k,m}}$ we may take $N =
\max N(\phi_{\alpha_{k,m}})$. A similar analysis as before shows
\eqref{eq:genh1Tn} also holds with $h_1$ replaced by $h_M$. The key
step in the induction is \bea & & \int_0^1 \left(h_M(y) - T_N
h_M(y)\right) g_{M+1}(x-y)\phi(x) dxdy \nonumber\\ & & \ \ = \
\int_0^1 \left(h_M(y) - T_N h_M(y)\right)
\sum_{k=1}^{r_{M+1}} w_{k,M+1} \phi(y+\alpha_{k,M+1})dy \nonumber\\
& & \ \ = \ \sum_{k=1}^{r_{M+1}} w_{k,M+1} \int_0^1 \left(h_M(y) -
T_N h_M(y)\right) \phi_{\alpha_{k,M+1}}(y)dy, \eea which, as the
$w_{k,M+1}$ sum to $1$, is less than $\gep/2$ in absolute value.
Arguing as in Theorem \ref{thm:mainL1} completes the proof.
\end{proof}

\appendix

%%%%%%%%%%%%%%%%%%%%%%%%%%%%%%%%%%%%%%%%%%%%%%%%%%%%%%%%%%%%%%%%%%%%%%%%%%%%%%%%%%
%%%%%%%%%%%%%%%%%%%%%%%%%%%%%%%%%%%%%%%%%%%%%%%%%%%%%%%%%%%%%%%%%%%%%%%%%%%%%%%%%%

\section{Comparison with Alternate Techniques}\label{sec:alternate}

We discuss an alternate proof of Theorem \ref{thm:mainbenfL1},
applying the standard Central Limit Theorem to the sum $\log_B |X_1|
+ \cdots + \log_B |X_M|$ and noting that as the variance of a
Gaussian increases to infinity, the Gaussian becomes uniformly
distributed modulo $1$. A significant drawback of a proof by the
Central Limit Theorem is the requirement (at a minimum) that the
variance of each $\log_B |X_m|$ be finite. This is a very weak
condition, and in fact many random variables $X$ with infinite
variance (such as Pareto or modified Cauchy distributions) do have
$\log_B |X|$ having finite variance; however, there are
distributions where $\log_B |X|$ has infinite variance.

To a density $f$ on $[0,\infty)$ we associate the density of the
mantissa, $f_B$. Explicitly, the probability that $X$ has first
digit (base $B$) in $[1,s)$ is just \be\label{eq:wraparound}
\int_1^s f_B(t)dt \ = \ \sum_{m=-\infty}^\infty \int_{1 \cdot B^m
\le x \le s \cdot B^m} f(x)dx. \ee

Let $X$ be the random variable with density \be\label{eq:falpha}
\twocase{f_\alpha(x) \ = \ }{\alpha / (x \log^{\alpha+1}
x)^{-1}}{{\rm if} $x \ge e$}{0}{{\rm otherwise.}} \ee This is a
probability distribution for $\alpha > 0$, and is a modification of
a Pareto distribution; see \cite{Mi} for some applications and
properties of this distribution. We study the distribution of the
digits base $e$; analogous results hold for other bases. The density
of $Y = \log X$ is $g(y) = \alpha y^{-(\alpha+1)}$ for $y \ge 1$ and
$0$ otherwise. For $\alpha \in (0,2]$ the random variable $Y$ has
infinite variance, and thus we cannot prove the Benford behavior of
products through the Central Limit Theorem; however, we can show the
random variable $X$ does satisfy the conditions of Theorem
\ref{thm:mainbenfL1}.

Let $F_{e,\alpha}$ be the cumulative distribution function of the
digits (base $e$) associated to the density $f_\alpha$ of
\eqref{eq:falpha}, and let $f_{e,\alpha}$ be the corresponding
density of $F_{e,\alpha}$. We assume $\alpha
> 1$ below to ensure convergence. By \eqref{eq:wraparound} we have
\be F_{e,\alpha}(s) \ = \ \int_1^s f_{e,\alpha}(t)dt \ = \
\sum_{m=0}^\infty \int_{1 \cdot e^m \le x \le s \cdot e^m}
f_\alpha(x)dx, \ee with $s \in [1,e)$. A simple integration gives
\be F_{e,\alpha}(s) \ = \ -\sum_{m=0}^\infty
\frac1{\log^{\alpha}(s\cdot e^m)}\ +\ \sum_{m=0}^\infty
\frac1{m^\alpha}; \ee note the second sum converges if $\alpha
> 1$. The derivative of the first infinite sum in the expansion of
$F_{e,\alpha}(s)$ is the sum of the derivatives of the individual
summands, which follows from the rapid decay of the summands (see,
for example, Corollary 7.3 of \cite{La}). Differentiating the
cumulative distribution function $F_{e,\alpha}$ gives the density
\be f_{e,\alpha}(s) \ = \ \alpha \sum_{m=0}^\infty \frac1{s
\log^{\alpha+1}(s \cdot e^m)}, \ \ \ s \in [1,e). \ee As $\alpha >
1$, for $m\neq 0$ the $m$\textsuperscript{th} summand is bounded by
$m^{-(\alpha+1)}$. Thus the series for $f_{e,\alpha}(s)$ converges
and is uniformly bounded for all $s$. A simple analysis shows that
the conditions of Theorem \ref{thm:mainbenfL1} are satisfied for
$\alpha \in (1,2]$.

The reason the Central Limit Theorem fails for densities such as
that in \eqref{eq:falpha} is that it tries to provide \emph{too
much} information. The Central Limit Theorem tries to give us the
limiting distribution of $\log_B |X_1 \cdots X_M| = \log_B|X_1| +
\cdots + \log_B |X_M|$; however, as we are only interested in the
distribution of the digits of $X_1 \cdots X_M$, this is more
information than we need.

%%%%%%%%%%%%%%%%%%%%%%%%%%%%%%%%%%%%%%%%%%%%%%%%%%%%%%%%%%%%%%%%%%%%%%%%%%%%%%%%%%
%%%%%%%%%%%%%%%%%%%%%%%%%%%%%%%%%%%%%%%%%%%%%%%%%%%%%%%%%%%%%%%%%%%%%%%%%%%%%%%%%%

%\ \\

\newpage

\end{document}